\documentclass[12pt]{article}

\usepackage{mathrsfs}
\usepackage{amsfonts}
\usepackage{amsmath, amsthm, amsfonts, amssymb}
\newtheorem{thm}{Theorem}[section]
\newtheorem{cor}[thm]{Corollary}
\newtheorem{lem}[thm]{Lemma}
\newtheorem{prop}[thm]{Proposition}

\newtheorem{Thm}{Theorem}[section]
\newtheorem{Cor}[Thm]{Corollary}
\newtheorem{Pro}[Thm]{Proposition}
\newtheorem{Lem}[Thm]{Lemma}

\theoremstyle{definition}

\numberwithin{equation}{section} \theoremstyle{remark}
\newtheorem{rem}{Remark}[section]
\def\<{\langle}
\def\>{\rangle}

\def\ra{\rightarrow}
\def\p{\partial}
\def\a{\alpha}
\def\g{\gamma}

\def\D{{\cal D}}
\def\O{\Omega}
\def\l{{\lambda}}

\def \sm{\setminus}

\def\-{\overline}

\def\e{\epsilon}

\def\a{\alpha}

\def\CC{\bf C}

\def\M*{\wt{M^*}}

\def\-{\overline}

\def\O{\Omega}

\def\D{\Delta}

\def\sm{\setminus}

\def\wt{\widetilde}
\def\ra{\rightarrow}

\def\a{\alpha}

\def\s{\sigma}
\def\D{\Delta}

\def\s{\mathcal}

\def\a{\alpha}

\def\beq{\begin{equation}}
\def\nneq{\end{equation}}
\def\beqn{\begin{eqnarray}}
\def\neqn{\end{eqnarray}}
\def\beqna{\begin{eqnarray*}}
\def\neqna{\end{eqnarray*}}
\def\bedis{\begin{displaymath}}
\def\nedis{\end{displaymath}}
\def\-{\overline}

\def\1#1{\overline{#1}}
\def\2#1{\widetilde{#1}}
\def\3#1{\widehat{#1}}
\def\4#1{\mathbb{#1}}
\def\5#1{\frak{#1}}
\def\6#1{{\mathcal{#1}}}

\def\C{{\4C}}
\def\CC{{\4C}}

\def\Label#1{\label{#1}}

\def\id{{\sf id}}

\def\Re{{\sf Re}\,}
\def\Im{{\sf Im}\,}

\def\phi{\varphi}
\def\eps{\varepsilon}

\def\a{\alpha}

\def\g{\gamma}
\def\G{\Gamma}
\def\D{\Delta}

\def\l{\lambda}

\theoremstyle{definition}\newtheorem{Def}[Thm]{Definition}

\theoremstyle{remark}
\newtheorem{Rem}[Thm]{Remark}
\newtheorem{Exa}[Thm]{Example}

\def\bl{\begin{Lem}}
\def\el{\end{Lem}}
\def\bp{\begin{Pro}}
\def\ep{\end{Pro}}
\def\bt{\begin{Thm}}
\def\et{\end{Thm}}
\def\bc{\begin{Cor}}
\def\ec{\end{Cor}}
\def\bd{\begin{Def}}
\def\ed{\end{Def}}
\def\br{\begin{Rem}}
\def\er{\end{Rem}}
\def\be{\begin{Exa}}
\def\ee{\end{Exa}}
\def\bpf{\begin{proof}}
\def\epf{\end{proof}}
\def\ben{\begin{enumerate}}
\def\een{\end{enumerate}}
\def\beq{\begin{equation}}
\def\eeq{\end{equation}}

\begin{document}

\
\bigskip\bigskip

\centerline{\bf  NON-EMBEDDABLE REAL ALGEBRAIC HYPERSURFACES}

\bigskip\bigskip
 \centerline{Xiaojun Huang\footnote{ Supported in part by
NSF-1101481} \ \  and\
Dmitri Zaitsev\footnote{Supported in part by the Science Foundation Ireland grant 10/RFP/MTH2878.}}


\begin{abstract}
We study various classes of real hypersurfaces
that are not embeddable into more special hypersurfaces in higher dimension,
such as spheres, real algebraic compact strongly pseudoconvex hypersurfaces or compact pseudoconvex hypersurfaces of finite type. We conclude by stating some open problems.
\end{abstract}

\bigskip\bigskip
\section{Introduction}
This paper is motivated by the following general problem:

\smallskip
{\em Given a real hypersurface $M$ in a complex manifold $X$, when can it be
(holomorphically) embedded into a more special real hypersurface $M'$ in a complex manifold $X'$ of possibly larger dimension? More specifically, which strongly pseudoconvex hypersurfaces can be embedded into a sphere?}
\smallskip

By a holomorphic map (resp.\ embedding) of $M$ into $M'$, we mean a
holomorphic map (resp.\ embedding) of an open neighborhood of $M$ in
$X$ into $X'$, sending $M$ into $M'$. In particular, it follows that
a hypersurface holomorphically embeddable into a sphere
$\4S^{2N-1}:=\{\sum_j |z_j|^2=1\}\subset\C^N$ is necessarily
strongly pseudoconvex and real-analytic. However, not every strongly
pseudoconvex real-analytic hypersurface can be even locally embedded
into a sphere, as was independently shown by Forstneric \cite{for}
and Faran \cite{fa}. These results, showing that such hypersurfaces
in general position are not embeddable into spheres, were more
recently further extended and strengthened by Forstneric \cite{For
2004} showing that they also do not admit transversal holomorphic
embeddings into a hyperquadric
$$\4H^{2N-1}_\ell:=
\{-\sum_{j\le\ell}|z_j|^2+\sum_{j>\ell} |z_j|^2=1\}\subset\C^N$$
of any signature $\ell$.
(By a transversal embedding $F$ we mean one
not sending the tangent space $T_pX$ into $T_{F(p)}\4H^{2N-1}_\ell$
for $p\in M$.)

Explicit examples
of non-embeddable strongly pseudoconvex
real-analytic hypersurfaces
were given by the second author \cite{Zo}
along with explicit invariants
serving as obstructions to embeddability.
In Theorem~\ref{ra} below
we give an example
of a {\em compact} strongly pseudoconvex
real-analytic hypersurface in $\C^2$
that does not admit any holomorphic embedding
into a sphere
(and more generally any transversal holomorphic
embedding into a hyperquadric).

The existence of non-embeddable real-analytic hypersurfaces suggests
to consider the embeddability problem for the more restricted class
of {\em real-algebraic} hypersurfaces, i.e.\ ones locally given by real
polynomial equations. In this line, Webster [We] showed in 1978 that any Levi-nondegenerate real-algebraic
hypersurface does in fact admit transversal holomorphic embeddings
into hyperquadrics of suitable dimension and signature. As a
consequence of the study of the Chern-Moser-Weyl tensor, Huang and Zhang
\cite{hz1} obtained  concrete  algebraic Levi non-degenerate
hypersurfaces with positive signature which can  not be holomorphic
embedded into a hyperquadric (with the same signature) of any
dimension.

During the Conference on Several Complex Variables and PDEs in Serra
Negra, Brazil, in August 2011, the authors observed that
  the strongly pseudoconvex (near $0$)
real-algebraic hypersurface  defined by
$$M:=\{z=(z_1,\cdots,z_n)\in {\C}^n: \Im{z_n}=\sum_{j=1}^{n-1}|z_j|^2-|z_1|^4\}, \quad n\ge 3,$$
is not locally (holomorphically)
embeddable into any sphere of any dimension nor into any closed
strongly pseudoconvex  real-algebraic hypersurface $M'\subset\C^N$
for any $N$. In fact, any such embedding would be algebraic by a
result of the first author \cite{h8} and hence would extend (as
holomorphic embedding into $M'$) to points of $M$ of mixed Levi
signature, which is impossible. In Theorem~\ref{1} below, we state a
generalization of this phenomenon leading to many simple examples of
strongly pseudoconvex real-algebraic hypersurfaces that are not
holomorphically embeddable even into any closed pseudoconvex
hypersurface $M'\subset\C^N$ of finite D'Angelo type.

One can  similarly construct the following locally non-embeddable
example in $\C^2$:
 $$M:=\{z\in {\C}^2: \Im{z_2}=|z|^2-|z|^4\},$$
where the proof is based on the observation that any potential
embedding would be extendable to ``large" Levi-degenerate sets,
which is impossible (see Theorem~\ref{2} below).
Along the same lines, we further study the property of a
class of real algebraic pseudoconvex hypersurfaces discovered by
Kohn and Nirenberg \cite{kn},
not to admit holomorphic supporting functions near certain
weakly pseudo-convex points (and hence not locally holomorphically
convexifiable near these points). We will prove a general
non-embeddability result in Theorem~\ref{KN-02} below, which is, in
addition to the Kohn-Nirenberg property, based on a property stated
in Proposition~\ref{KN-03}, which roughly says that in certain
situations a holomorphic extension of a local embedding from $M$
into $M'$ even along paths outside $M$ still sends $M$ into $M'$.
Proposition~\ref{KN-03} is a generalization of what is called the
invariant property for holomorphic corespondences in the literature
(see \cite{hu3}). However, our proof here is   more geometric and
also   simpler even in the case considered in \cite{hu3}. Our
general non-embeddability theorem immediately leads to many examples
of compact pseudoconvex real-algebraic hypersurfaces, strongly
pseudoconvex away from a single point, that are not locally
holomorphically embeddable into any compact strongly pseudoconvex
real-algebraic hypersurface of any dimension. We also mention
recent related preprint by Ebenfelt and Son [ES].


\medskip
We next address  the related problem for hypersurfaces of positive (mixed)
Levi signature. That is, whether there exists a {\em compact Levi
nondegenerate real-algebraic hypersurface of signature $\ell>0$}
that is not transversally embeddable into a hyperquadric
$\4H^N_\ell$ of higher dimension but the same signature $\ell$. Note
that Webster's result \cite{We1} shows that without the signature
restriction, such an embedding is always possible. However, based on
a monotonicity property of the {\em Chern-Moser-Weyl tensor}
 \cite{hz1}
and algebraicity results in \cite{Hu2} (see also \cite{Z} and
\cite{CMS}), we give in \S\ref{mixed} below examples of compact
real-algebraic Levi-nondegenerate
 hypersurfaces of positive Levi signature in the projective space
that are not transversally locally holomorphically embeddable
into any hyperquadric of any dimension but the same signature.
(Note that there is no compact hypersurface in $\C^n$ with positive signature,
since any such hypersurface must have a strongly pseudoconvex point.)

Finally we mention some open problems in the last section.

\section{Hypersurfaces not embeddable
into certain  real-algebraic hypersurfaces}

We first recall that a smooth real hypersurface in an open subset
$U$ of $\C^n$ is called real algebraic, if it has  a real-valued
polynomial defining function. A real algebraic hypersurface in $U$
has an extension to a real analytic variety in $\C^n$, which may
possess singularities. Of course,   all real algebraic hypersurfaces
in $\C^n$ are automatically smooth and closed.

\subsection{A compact strongly pseudoconvex real-analytic
hypersurface not embeddable into any strongly pseudoconvex
real-algebraic hypersurface} In \cite[Corollary 1.2]{Z}, the second
author gave an explicit example of a germ of real-analytic strongly
pseudoconvex hypersurface in $\C^2$ that is not transversally
holomorphically embeddable into any Levi-nondegenerate
real-algebraic hypersurface. By following verbatim the proof of [Z2,
Corollary 1.2], one has:

\bt\label{ra}
For $0<\eps<<1$,
$$
M=\left\{(z,w)\in\C^2 : |z|<1, \, |w|^2 +  |z|^2 + \eps\Re
\sum_{k\ge 2} z^k \bar z^{(k+2)!} = 1/2 \right\}
$$
is a compact strongly pseudoconvex hypersurface that does not
admit a nontrivial holomorphic embedding into any Levi-nondegenerate
real-algebraic hypersurface in $\C^N$. In particular, $M$ is not
holomorphically embeddable into any strongly pseudoconvex
real-algebraic hypersurface in $\C^N$. \et

\subsection{Strongly pseudoconvex real-algebraic hypersurfaces not
embeddable into  finite type real-algebraic hypersurfaces}

Recall that a real-analytic hypersurface
is of {\em finite D'Angelo type}
if and only if it does not contain
any complex curve.
By a point of {\em mixed Levi signature}
we mean a point of a real hypersurface,
where the Levi form (for a choice of conormal) has both positive and negative eigenvalues.
We next state the following:
 \bt\Label{1} Let $M\subset\C^{n+1}$ be
a connected real-algebraic hypersurface with a point of mixed Levi
signature. Then any holomorphic map sending an open subset of $M$ into any
closed pseudoconvex finite D'Angelo type real-algebraic
hypersurface $M'\subset\C^{N+1}$ is constant. \et

\bpf Obviously, $M$ has nonzero Levi form on a dense subset. By
 [Z1], [CMS] (or by [Hu1] when the target is strongly
pseudoconvex), any holomorphic map $F$ sending an open subset of $M$
into $M'$ is complex-algebraic. Since $M$ is connected and since the
branching variety of $F$ is of complex codimension  one (if $F$ is
not a single-valued), we can extend $F$ along a path to a
neighborhood of a point $p\in M$ of mixed Levi signature, still
sending $M$ into $M'$. Since $M'$ is pseudoconvex, $F$ must be
constant near $p$ and hence it is constant. \epf

\begin{Exa} Consider the hypersurface
\begin{equation}
\label{00-0}
 M:=\{z\in {\bf C}^n: n\ge 3,\  \Im{z_n}=|z|^2-|z_1|^4\}.
\end{equation}
Then no open piece of $M$ can be holomorphically embedded into any closed
real-algebraic  hypersurface of finite D'Angelo type in $\C^N$ for any $N$.
\end{Exa}

\subsection{Hypersurfaces with large Levi-degenerate set}

\bt\Label{2} Let $M\subset\C^{n+1}$ be a connected real-algebraic
hypersurface, Levi-nondegenerate at some point, whose set of
Levi-degenerate points contains a real-analytic submanifold that is
generic in $\C^{n+1}$. Then any holomorphic map sending an open
subset of $M$ into any  strongly pseudoconvex real-algebraic
hypersurface $M'\subset\C^{N+1}$ is constant. \et

\bpf It follows from the assumptions that the set $S\subset M$ of Levi-degenerate
points is a generic real-analytic submanifold near some point $p\in
S$. As in the proof of Theorem~\ref{1} any holomorphic map $F$ sending an open
subset of $M$ into $M'$ extends holomorphically and algebraically
into an open neighborhood of a point $p\in S$ as above, still
sending $M$ into $M'$. (Note that algebraicity here already follows
from \cite{h8}.) Since $M'$ is strongly pseudoconvex, the extension
$F$ must have rank less than $n+1$ for all $q\in S$ near $p$. Since
$S$ is a generic submanifold of $\C^{n+1}$, the rank of $F$ is less
than $n+1$ in an open neighborhood of $p$. On the other hand, $F$ is
either constant or of full rank $n+1$ at any Levi-nondegenerate
point of $M$. Hence $F$ must be constant. \epf

\be A simple example of $M$ satisfying the assumptions of
Theorem~\ref{2} is
$$M:=\{(z,w)\in \C^2 : \Im w = |z|^2 - |z|^4\},$$
which is strongly pseudoconvex near $0$. The Levi-degenerate set
here is $\{(z,w)\in M: |z|=1/2\}$ and hence Theorem~\ref{2} applies.
\ee

\be An example of a compact pseudoconvex $M$ satisfying the
assumptions of Theorem~\ref{2} is the following boundary of a
Reinhardt domain:
$$M:=\{ (z,w)\in\C^2 : (|z|^2+|w|^2)^4 +(|z|^2-|w|^2)^4 =1\}.$$
Away from the cross $zw=0$, $M$ is locally biholomorphically
equivalent to
$$\2M:=\{ (z,w)\in\C^2 : (|z|+|w|)^4 +(|z|-|w|)^4 =1\},$$
under the finite holomorphic map $(z,w)\mapsto (z^2,w^2)$. The real
part of $\2M$ is the rotated convex curve $\{x^4+y^4= 1/4\}$
whose real part is convex. Then $M$ is pseudoconvex and is
Levi-degenerate along the generic submanifold $\{(z,w)\in M :
|z|=|w|\}$. \ee

\section{Non-embeddable Kohn-Nirenberg type domains}\label{KN}

\be Given integers $0<l<k$, consider the following famous
Kohn-Nirenberg  domain \cite{kn}. (In the  paper of Kohn-Nirenberg
\cite{kn}, though only the case with $l=1,k=4, c=\frac{15}{7}$ was
studied, the result in their paper holds, with the same argument,
for the following more general domain which we still call the
Kohn-Nirenberg domain):
\begin{equation}\label{31}
\Omega=\{(z,w)\in {\bf C}^2: \rho=-\Im{w}+z^k\bar z^k+c\Re(z^{l}\bar
z^{2k-l})<0\}, \quad 2<|c|<\frac{k^2}{l(2k-l)}.
\end{equation}
 Also notice that the boundary of $\Omega$ is of
type $2k$ at $0$ and of bi-type $(l,2k-l)$. It is easily seen that
$\O$ is smooth. Since the Levi form of $\p\Omega$ is positive over
$\p \O\sm L_0$ with $L_0:=\{\Im{w}=0, z=0\}$, and is semi-definite
along $L_0$, we see that $\O$ is strongly pseudoconvex away from
$L_0$ and is weakly pseudoconvex of finite type along $L_0$.
Kohn and Nirenberg \cite{kn} proved the
following basic feature of the boundary of $\Omega$ that we call here
{\em Kohn-Nirenberg property}:
\ee

\begin{Def}\label{k-n}
A real hypersurface $M\subset\C^n$ is said to satisfy the {\em
Kohn-Nirenberg property} at a point $p\in M$ if  for any holomorphic
function $h\not\equiv0$ in any neighborhood $U$ of $p$ in $\C^n$ with
$h(p)=0$, the zero set $\6Z$ of $h$ intersects $M$ transversally at
some smooth point of $\6Z$ near $p$.
\end{Def}

In particular, a hypersurface with the Kohn-Nirenberg property at a
point is always minimal at that point. (We mention that when $M\cap
{\cal Z}$ separates ${\cal Z}$, it has Hausdorff codimension one and
thus must be generically smooth in ${\cal Z}$. ) See also
Example~\ref{compact} for {\em compact} hypersurfaces with the
Kohn-Nirenberg property. The argument in \cite{kn} is very general
and can be used to obtain further classes of hypersurfaces
satisfying the Kohn-Nirenberg property. We mention the paper by M.
Kolar \cite{ko} for a discussion of the similar but different
property of local holomorphic convexifiability.

We shall also consider local
holomorphic supporting functions:

\bd
A subset $M'\subset\CC^N$ is said to admit
{\em local holomorphic supporting functions}
if for each $q\in M'$,
there is a neighborhood $\O$ of $q$ in $\CC^N$ and a holomorphic function
$h$ in $\O$ such that $h(q)=0$ but $\Im {h(z)}<0$
for $z\in M'\cap\O$, $z\ne q$.
\ed

\br
 In particular,  when $M'$ is a smooth hypersurface of finite D'Angelo type and
  locally holomorphically convexifiable, it admits local
holomorphic supporting functions. This is a  consequence of a result
of McNeal on the equivalence of linear type and D'Angelo type for
convex domains. (See [DF2], for instance).


Theorem  \ref{KN-02} below implies that no open piece of the
boundary of the classical Kohn-Nirenberg domain can be mapped by a
non-constant holomorphic map into any connected compact smooth
algebraic hypersurface in $\C^n$,
that is locally holomorphically convexfiable.
\er

\bigskip
\be\Label{compact}
Consider the following  compactified Kohn-Nirenberg type domain:
\begin{equation}\label {0101}
\Omega=\{(z,w)\in {\C}^2: \rho=\epsilon
(|z(w-1)|^2+|z|^{2k}+c|z|^{2l}\Re{z^{2k-2l}})+|w|^2+|z|^{2k+2}-1<0\}.
\end{equation}
where $0<\e<<1$ and  $l,k, c$ as in \eqref{31}.
 Then $\O$ is a smoothly
bounded real-algebraic domain, which is pseudoconvex and strongly
pseudoconvex away from $p_0:=(0,1)$. Since the principal terms in
$\rho$ at $p_0$ are the same as those in the classical
Kohn-Nirenberg domain case, one still has the Kohn-Nirenberg
property at $p_0$ by the same argument.
Again, by Theorem \ref{KN-02} below, no open piece of $\p \O$ can be
mapped by a non-constant holomorphic map into a
smooth compact algebraic hypersurface $M'$, that admits  local
holomorphic supporting functions.

To get a similar  higher dimensional example with the Kohn-Nirenberg
property, we need only to find  one which includes the boundary of
the domain in (\ref{0101}) as its CR submanifold. For instance, the
boundary of the following domain  serves this purpose:
\begin{multline}\label {010101}
\{(z,w)=(z_1,z',w)\in {\C}\times{\C}^{n-2}\times{\C}:\\
 \rho=\epsilon
(|z_1(w-1)|^2+|z_1|^{2k}+c|z_1|^{2l}\Re{z_1^{2k-2l}})+|w|^2+
|z'|^2+|z_1|^{2k+2}<1\}.
\end{multline}

 \ee

\bt\label{KN-02} Let $M\subset {\CC}^n$ $(n>1)$ be a connected
minimal real-algebraic hypersurface, which has the Kohn-Nirenberg
property at some point. Let $M'\subset\CC^N$ be a compact
real-algebraic subset admitting local holomorphic supporting
functions at each point. Then any holomorphic map sending an open
piece of $M$ into $M'$
 is constant.
 \et

The proof is broken up in a sequence of lemmas.

\begin{lem}\label{fundamental}
Let $U\subset\C^n$ be a simply connected open set and $\6S\subset U$
a closed complex analytic subset of codimension one. Then for $p\in
U\setminus \6S$,
 the fundamental group $\pi_1(U\setminus\6S,p)$
is generated by loops obtained by concatenating paths
$\g_1$, $\g_2$, $\g_3$, where $\g_1$ connects $p$ with a point
arbitrarily close to a smooth point $q_0\in \6S$, $\g_2$ is a loop
around $\6S$ near $q_0$ and $\g_3$ is $\g_1$ reversed.
\end{lem}
Here, by saying that $\g_2$ goes around $q_0$, we mean there is a
closed  embedded real 2-disk $\overline{D}$ in $U$ such that $\g_2$
is the boundary of $\overline{D}$ and $D$ intersects ${\cal S}$ only
and transversally  at $p_0$.

 \bpf Replacing $U$ by $U\sm Sing({\cal S})$ if needed, we can
assume that $S$ is smooth. Here $Sing({\cal S})$ is the
singular set of ${\cal S}$, which  has codimension at least two, hence $U\sm
Sing({\cal S})$ is still simply-connected. Take any loop $\g\in \pi_1(U\setminus\6S,p)$.
Since $U$ is simply connected, $\g$ can be contracted to $p$ inside
$U$, i.e.\ $\g$ viewed as a map from $S^1:=\{|z|=1\}\subset\C$ into
$U$ can be continuously extended to the disk $\1\D:=\{|z|\le 1\}$.
Using Thom's transversality, the disk extension can be approximated
by a smooth immersion $\G\colon \1\Delta\to U$ such that $\G|_{S^1}$
is a smooth Jordan loop defining the same class in
$\pi_1(U\setminus\6S,p)$ as $\g$, and such that $\G(\Delta)$
intersects $\6S$ transversally at finitely many smooth points.
Since the fundamental group of the disk $\Delta$
minus finitely many points is generated by loops
going around single points, it is easy to see
that $\G|_{S^1}$ and hence $\g$ is generated
by loops inside $\G(\1\Delta)$ as described in the lemma.
\epf

Let $M\subset U(\subset\C^n)$ be a closed real-analytic subset
defined by a family of real-valued real analytic functions
$\{\rho_\a(z,\-{z})\}$. Assume that the complexification
$\rho_\a(z,\xi)$ of $\rho_\a(z,\-{z})$ is holomorphic over $U\times
\hbox{conj}(U)$ with
$$\hbox{conj}(U):=\{z:\-{z}\in U\}$$ for each
$\a$. Then the complexification $\6M$ of $M$
is  the complex-analytic subset in $U\times\hbox{conj}(U)$ defined
by $\rho_\a(z,\xi)=0$ for all $\a$.
Then for $w\in\C^n$, the Segre variety of $M$ associated with the
point $w$  is defined by $Q_w:=\{z : (z,\bar w)\in\6M\}$. Recall the
basic properties of the Segre varities: $z\in Q_w \iff w\in Q_z$ and
$z\in Q_z \iff z\in M$. (See \cite{hu3} for more related notations
and definitions.)

\begin{lem}\label{u-choice}
Let $M\subset\CC^n$ be a minimal real-analytic hypersurface at a
point $p_0\in M$. Then there exist small open neighborhoods $U, \2U$
of $p_0$ in $\CC^n$ with  $U\subset\subset \2U$  such that the
following holds:

\begin{enumerate}
\item For every $z\in U$, the Segre variety $Q_z$
is a nonempty closed connected smooth hypersurface in $\2U$.
\item There is no complex hypersurface $H\subset U$
such that $Q_z \equiv Q_w$, when restricted to $\wt{U}$,  for all
$z,w\in H$.
\end{enumerate}
\end{lem}

\bpf Let $M$  be a real analytic hypersurface near $p_0$ as in the
lemma with a real analytic defining function $\rho$ near $p_0$. (1)
is a direct consequence of the implicit function theorem and is standard in the literature.

We prove (2) by contradiction assuming there exists $H$ as in the lemma.
Suppose $p_0\in H$. Since $p_0\in Q_{p_0}$, and
 for any $w\in H$,
we must have $Q_w\equiv Q_{p_0}$. Hence $p_0\in Q_w$ and therefore
$w\in Q_{p_0}\equiv Q_w$. Hence $w\in Q_w$ and thus $H\subset M$,
contradicting nonminimality of $M$.

For $H$ general, and for $z\in H$  $q\in Q_z$, we have $q\in
Q_w\equiv Q_z$ for any $w\in H$. Therefore $w\in Q_q$, and thus
$H\subset Q_q$, which gives that $H=Q_q$. Hence, by the property of
$H$, we see that $E_q:=\cup_{z\in Q_{q}} Q_z=Q_{z_0}$ for any
$z_0\in H$ and thus is a complex hypersurface.

On the other hand, assume without loss of generality that $p_0=0$
and $\frac{\p{\rho}}{\p z_n}(0)\not = 0$ for a real-analytic
defining function of $M$. Then there is a holomorphic function
$\Psi$ in its variables such that $E_q$ is defined, in the
$(\xi',\xi_n)$-coordinates, near $0$ by $\xi_n=\Psi(\-{z'},q,\xi')$
with parameter $z'\approx 0$. The latter notation here means that
$z'$ is sufficiently close to $0$ and we shall use it in the sequel.
Now, suppose that the statement in (2) fails no matter how we shrink
$U$. Then we have a sequence $q\ra 0$ such that  $E_q$ is simply
defined by $\xi_n=\Psi(0,q,\xi')$. Passing to the limit, we get
$E_0$ is defined by $\xi_n=\Psi(0,0,\xi')$. This contradicts the
minimality as argued above.
\epf

\begin{lem}\label{H}
Let $M\subset\CC^n$ be a minimal real-analytic hypersurface at a
point $p_0\in M$ and $\6S$ a closed proper  complex analytic subset
 in a neighborhood of $p_0$ with $p_0\in {\cal S}$. Then there exists a small (simply-connected) open
neighborhood $U$ of $p_0$ in $\CC^n$, such that the following holds.

Take $p\in (M\cap U)\setminus\6S$ and let
$\g\in\pi_1(U\setminus\6S,p)$ be obtained by concatenation of
$\g_1,\g_2,\g_3$ as in Lemma~\ref{fundamental}, where $\g_2$ is a
small loop around $\6S$ near a smooth point $q_0\in \6S\cap U$. Then
$\g$ can be slightly perturbed to a homotopic loop $\2\g(t)$ in
$\pi_1(U\setminus\6S,p)$ such that there exists a null-homotopic
loop $\l(t)$ in $\pi_1(U\setminus\6S,p)$ with
$(\l(t),\overline{\2\g(t)})$ contained in the complexification $\6M$
of $M$ for all $t$. Also, for any element $\hat{\g}\in
\pi_1(U\setminus\6S,p)$, after a small perturbation of $\hat{\g}$ if
needed,
we can find a null-homotopic loop  $\hat{\l}\in
\pi_1(U\setminus\6S,p)$ such that $(\hat{\g},\-{\hat{\l}})\subset
\6M$.
\end{lem}

\bpf Let $p_0\in U\subset\2U$ be satisfying the conclusion of
Lemma~\ref{u-choice}. Shrinking $U$ if necessary, we may assume that
there exists a real analytic (reflection) map $\6R\colon U\to {U}$
with ${\cal R}^2=\id$, $\6R|_M=\id$ and $({\6R(z)},\1{z})\in \6M$
for all $z\in U$, namely, ${\cal R}(z)\in Q_z$. In fact, the map
$\6R$ can be obtained by slicing $M$ transversally by a family of
parallel complex lines $\{L\}$ near $p_0$ and then taking the
Schwarz reflection about $M\cap L$  inside each $L$ of the family.
More precisely, let $L_{0}$ be a complex line through $p_{0}$
intersecting $M$ transversally at $p_{0}$. Then sufficiently small
neighborhood $U$ of $p_{0}$ is foliated by lines $L$ parallell to
$L_{0}$, which still intersect $M$ transversally. Shrinking $U$
suitably, we may assume that the Schwarz reflection about $M\cap L$
is defined in $U\cap L$ and leaves the latter invariant. Then define
$\6R$ to be the Schwarz reflection along each line $L$. (We can of
course arrange $U$ such that for any line $L$, the pair of the
reflecting points with respect to $U\cap L$ stays inside $U$ and
thus ${\cal R}(U)=U.$)


We now claim that we can slightly perturb $q_0\in\6S$ and
 the direction of the parallel lines (and hence $\6R$)
such that $\6R(q_0)\notin \6S$. Indeed, by  (2) of
Lemma~\ref{u-choice} applied to  $H=\6S$, we conclude that
$Q_{q}\not\equiv Q_{q'}$ for two generic $q,q'\in \6S$ arbitrarily
close to $q_0$. Then either $Q_{q_0}$ contains points away from
$\6S$ arbitrarily close to $q'_0:=\6R(q_0)$ or an open piece of
$Q_{q_0}$ is contained in $\6S$. But the latter case together with
$Q_q\not\equiv Q_{q'}$ with $q,q'(\in {\cal S})\approx q_0$ implies
that $Q_q$ cannot contain an open piece of $\6S$ for a generic
$q\approx q_0$. Then we can choose such $q$ and $q'\in
Q_q\setminus\6S$ arbitrarily close to $q_0$ and $q_0'={\cal
R}(q_0)$, respectively. Considering the line through $q$ and $q'$
and using the  lines parallel to this one to redefine ${\cal R}$,
this proves the claim.

After slightly perturbing $q_0$ and $\6R$ as in the above, it
follows that there exists a sufficiently small open ball $\Omega$
containing $q_0$ such that $\6R(\Omega)\cap \6S=\emptyset$. Then the
paths $\g_1, \g_2, \g_3$ can be perturbed homotopically into
$\2\g_1,\2\g_2,\2\g_3$ respectively, where $\2\g_1$ connects $p$
with a point in $\Omega$, $\2\g_2$ is a loop around $\6S$ inside
$\Omega$ and also sufficiently close to $q_0$, and   $\2\g_3$ is
$\2\g_1$ reversed such that the loop $\2\g$ obtained by
concatenation of $\2\g_1,\2\g_2,\2\g_3$ is homotopic to $\g$ in
$\pi_1(U\setminus\6S,p)$ and we can take $\l(t):=\6R(\wt\g(t))$. (Of
course, we may need to slightly perturb $\wt\g_1$ to make sure that
$\l$ avoids $\6S$.)  Then $\l$ is null-homotopic in
$\pi_1(U\setminus\6S,p)$ since $\6R(\Omega)$ does not intersect
$\6S$.

The last statement in the lemma follows from the symmetric property
of the reflection map (Segre varieties) and
what we just proved.
 The proof is complete. \epf

We now choose $\6R$ as in the above proof above, defined in a
neighborhood of a point $p_0\in M$.

\begin{prop}\label{KN-03}
Let $\O\subset\CC^n$, $V\subset\CC^N$ be connected open sets,
$M\subset \O$ a real analytic hypersurface, $M'\subset V$  a
real-analytic subset defined by a set of real valued  real analytic
functions $\{\rho_\a\}$ over $V$, ${\6S}\subset \O$ a
 proper closed complex analytic subset
 and ${\cal F}\subset (\O\setminus {\6S})\times V$ a
complex submanifold whose projection to $\O\setminus {\6S}$ is a
finite sheeted covering. Suppose that:

\begin{enumerate}
\item $M$ is minimal at $p_0\in M$.
\item The complexification $\rho_\a(z,\xi)$  for each $\a$ is holomorphic over $V\times \rm{conj}(V)$.
\end{enumerate}
Then there exists a neighborhood $U$ of $p_0$ in $\O$,
depending only on $M$ and $p_0$, such that if a certain local branch
$F$ of $\6F$, defined over a subdomian $D\subset U\sm {\cal S}$ with
$D\cap M\not =\emptyset$ sends a $D\cap M$ into $M'$, then any other
branch of $\6F$ obtained by continuing $F$ along paths in
$U\setminus\6S$ also sends $M$ into $M'$. Equivalently, if ${\cal
F}'$ is the connected component of ${\cal F}\cap(U\setminus {\cal S})\times V$ containing the graph
of $F$, then for any $(z,w)\in {\cal F}'$ with $z\in M$, we have $w\in M'$.
 More generally, slightly perturbing $\6R$  if needed, we have
$F_1(Q_z\cap O({\cal R}(z)))\in Q'_{F_2(z)}$ for any $(z,F_1(z)),(z,F_2(z))\in {\cal F}'$ with $z, {\cal R}(z)\in U\sm {\cal S}$.
Here for any $w\in V$, $Q'_w:=\{z\in V: \rho_\a(z,\-{w})=0\ \forall
\a\}$, and $O(a)$ denotes a small neighborhood of $a$ in ${\CC}^n$.

\end{prop}

\bpf  Let $U$  be a simply connected neighborhood of $p_0$ in $\Omega$.
We need
only consider the case when ${\cal S}$ is of codimension one in
$\Omega$, for, otherwise, $U\sm {\cal S}$ is also simply connected.
Hence, the continuation of $F$ along curves in $U\sm {\cal S}$
defines a holomorphic map over $U\sm {\cal S}$ and all the
properties stated in the Proposition follows easily. We also choose
$U\subset\O$ such that the conclusion of Lemma~\ref{H} is satisfied.
In addition, we can choose $U$  such that $M\cap U$ is connected and
minimal.

 Denote by $F\colon
D\subset U\setminus \6S\to V$ a local branch of $\6F$,
where $D$ is a domain with $D\cap M\not =\emptyset$,
sending $D\cap M$ into $M'$. Let  $p\in D\cap M$.
Then $(F,\1F)\colon D\times \hbox{conj}(D) \to
V\times\hbox{conj}(V)$ sends an open neighborhood of $(p,\bar p)$ in the complexificaiton $\6M$ of $M$ into $\6M'$.

Let $F_1\colon D_1\subset U\setminus\6S\to V$ be a branch of $\6F$
with some point $p_1\in M\cap D_1$, obtained by continuing $F$ along a path in
$U\setminus\6S$, connecting $p_1$ with $p$. Since $M\cap U$ is
connected and minimal, $(M\cap U)\setminus\6S$ is also connected.
Hence there exists a path $\g$ in $(M\cap U)\setminus\6S$ connecting
$p$ with $p_1$. Then by the analyticity of $M'$, the branch $F_2$ of
$\6F$ obtained by continuing $F$ along $\g$ is sending a
neighborhood of $p$ in $M$ into $M'$.

Hence $(F_2,\-{F_2}):=(F_2(\cdot),\1{F_2(\-\cdot)})$ sends a
neighborhood of $(p_1,\1{p}_1)$ in $\6M$ into $\6M'$. Now the branch
$F_1$  is obtained from $F_2$ by continuation along a certain loop
$\g$ in $\pi_1(U\setminus\6S,p_1)$. Notice that ${\cal R}^2=\id$. By
Lemma \ref{fundamental} and Lemma \ref{H},  slightly perturbing $\g$
and $\cal R$  if needed, we can assume that $\l(t)={\cal R}(\g(t))$
is a null homotopic loop in $U\sm {\cal S}$. Notice that $\g={\cal
R}(\l)$. Applying the holomorphic continuation along the loop
$(\g,\-{\l})$ in $\6M$ for $\rho_\a(F_1,\1F_1)$ for each $\a$, one
concludes by the uniqueness of analytic functions that
$(F_1,\1F_2)$ also sends a neighborhood
of $(p_1,\1{p}_1)$ in $\6M$ into
$\6M'$.  Namely, for any $z$ near $p_1$, we have
$F_1(Q_z\cap O(\6R(z)))\subset Q'_{F_2(z)}$,
where $O(a)$ as before denotes a small neighborhood of $a$.
Now, applying
the holomorphic continuation  along the loop $(\l, \-{\g})$ in $\6M$
for $\rho_\a(F_1,\1{F}_2)$ for each $\a$, one concludes by the
uniqueness of analytic functions that
$(F_1,\1F_1)$  sends a neighborhood of $(p_1,\1{p}_1)$ in $\6M$ into
$\6M'$. In particular, $F_1$ maps a neighborhood of $p_1$ in $M$ into
$M'$. (Cf.\ Lemma 2.1 of \cite{hj4}).
The last assertion in the proposition  can be proved with a similar
argument based on the holomorphic continuation and the uniqueness
property for analytic functions.\epf

\begin{lem}\label{KN-04}
For an open set $U\subset \CC^n$, consider the complex analytic
subset \beq\label{404'} \6F:= \left\{ (z,w)\in U\times\CC :
\sum_{l=0}^{m}a_{l}(z)w^l = 0 \right\},
 \nneq
where $a_0(z),\ldots,a_m(z)$ are holomorphic functions in $U$ that
do not simultaneously vanish on a (possibly singular)  complex
hypersurface. Suppose that $M\subset U$ is  a real-analytic
hypersurface and $C>0$ is a constant such that $|F(z)|\le C$ for any
branch $F\colon\Omega\to \CC$ of $\6F$ and any $z\in M\cap \Omega$.
Write $\6S':=\{z\in U: a_{m}(z)=0\}$. Then $M\cap \6S'$ is contained
in a complex analytic subset of $\6S'$ of positive codimension.
\end{lem}


\bpf Since $a_0(z),\ldots,a_m(z)$  do not simultaneously vanish on a
(possibly singular) complex hypersurface, for each non-empty
irreducible component  $\6C$  of $\6S'$, there exists $j<m$ such
that $a_j(z)$, does not vanish identically on $\6C$. Hence
$\{a_j=0\}$ defines a complex analytic subset of $\6S'$ of positive
codimension.

We claim that $M\cap \6C\subset \{a_j=0\}$.
Indeed, suppose on the contrary, there exists $z_0\in M\cap\6C$ with $a_j(z_0)\ne 0$.
Since $M$ is a real hypersurface, there exists
a sequence $z_k\in M\sm \6S'$ converging to $z_0$ as $k\to\infty$
such that $\6F$ has $m$ branches (counted with multiplicity)
around each $z_k$.
Since all branches of $\6F$
are uniformly bounded on $M$ by our assumption,
the same is true for their symmetric functions.
In particular, $\frac{a_j(z)}{a_m(z)}$ is also bounded
along $z_k$.
On the other hand $a_m(z_0)=0$, $a_j(z_0)\ne0$
imply $\frac{a_j(z_k)}{a_m(z_k)}\to \infty$ as $k\to\infty$,
which is a contradiction.

This proves the claim and therefore,
$M\cap\6C$ is contained
in the set $\{a_j=0\}$ of positive codimension.
Since $\6C$ is an arbitrary irreducible
component of $\6S'$, the proof is complete.
\epf

\begin{cor}\label{kn-cor}
In addition to the assumptions of Lemma~\ref{KN-04},
suppose that $M$ satisfies the
Kohn-Nirenberg property at $p\in M$
(see Definition~\ref{k-n}).
Then $a_m(p)\ne 0$.
\end{cor}

\bpf
Indeed, otherwise by the Kohn-Nirenberg property,
the zero set $\6S'$ of $a_m(z)$
must intersect $M$ transversally
at some smooth point.
The latter implies that $M\cap \6S'$
is a real hypersurface in $\6S'$ at such point,
which contradicts the conclusion of Lemma~\ref{KN-04}.
\epf

%

\bigskip

\bpf[Proof of Theorem \ref{KN-02}]
Assume the hypotheses in Theorem \ref{KN-02}.
Suppose that there is a non-constant holomorphic map $F$ sending
an open piece $M$ into $M'$.
By a  result of Diederich-Fornaess \cite{DF1}, $M'$ does not
contain non-trivial holomorphic curves.
Since $M$ is minimal, by
\cite{Z} and \cite{CMS}, $F$ is complex algebraic. In particular,
$F$ extends holomorphically along any path away from a proper
complex algebraic subset $\6S\subset\CC^n$. We need only prove the
theorem assuming that ${\cal S}$ is a codimension one complex
analytic variety near $p_0\in M$ with the Kohn-Nirenberg property.

Since $M$ is minimal and connected, $M\sm \6S$ is also connected.
Then $F$ has holomorphic extensions to points of $M\sm\6S$
arbitrarily close to a point $p_0\in M$ with the Kohn-Nirenberg
property, sending $M$ into $M'$. Now Proposition~\ref{KN-03} implies
that there exists a neighborhood $U$ of $p_0$ in $\CC^n$ and an
extension $\2F$ of $F$ to a point in $M\cap U$ such that any
extension of $\2F$ along a path in $U\sm\6S$ sends $M$ into $M'$.
Consider the ($n$-dimensional) Zariski closure
$\2{\6F}\subset\CC^n\times\CC^N$ of the graph of $\2F$ and denote by
$\3{\6F}$ the analytic irreducible component of $\2{\6F}\cap
(U\times\CC^N)$ containing the graph of $\2F$. In particular,
$\3{\6F}\sm (\6S\times\CC^N)$ is connected and therefore each branch
of $\3{\6F}$ away from $\6S$ sends $M$ into $M'$.

Since $M'$ is compact, it follows that all branches of $\3{\6F}$ are
uniformly bounded on $M$. Then Corollary~\ref{kn-cor} implies that,
after possible shrinking $U$ around $p_0$, $\3{\6F}$ becomes
bounded. Furthermore, by further shrinking $U$, we may assume that
$\3{\6F}\cap (\{p_0\}\times\CC^N) = \{(p_0,p'_0)\}$ for some
$p'_0\in M'$.
%

Since $M'$ has local holomorphic supporting functions,
there exists a holomorphic function $h$ in a neighborhood of $p'_0$ in $\CC^N$
such that $h(p'_0)=0$ and $\Im h<0$ on $M'$ away from $p'_0$.
Let $F_1,\ldots, F_m$ be local branches of $\3{\6F}$
at $z\in U\sm \6S$, counted with multiplicity.
Define
$h^*:=\sum_{j=1}^{m}h\circ F_j$.
Then $h^*$ is well-defined away from $\6S$
and  extends
holomorphically to $p_0$
with $h^*(p_0)=0$.
Furthermore, since all branches of $\3{\6F}$ send $M$ into $M'$,
we have
 $\Im{h^*(z)}< 0$ for generic  $z\in M$
 unless $\3{\6F}\cap (\{z\}\times\CC^N) = \{(z,p'_0)\}$.
 Since $F$ is assumed to be non-constant, so is $\2F$.
 Hence there exist points $z\in M$ arbitrarily close to $p_0$
 with $\Im h^*(z)<0$.
 In particular, $h^*\not\equiv0$ and hence, by the
 Kohn-Nirenberg property,
 the zero set $\6Z:=\{h^*=0\}$
 intersects $M$ transversally at some smooth points of $\6Z$
 arbitrarily close to $p_0$.

 Since $M$ is minimal,
 one-sided neighborhood $D$ of $p_0$
 is filled by small analytic disks in $U$ attached to $M$
 by a result of Tr\'epreau \cite{Tr} (see also \cite{Tu}).
Therefore we have $\Im h^*\le 0$ in $D$ by the maximum principle.
Since $\6Z$ intersects $M$ transversally at some points close to
$p_0$, it also intersects $D$. That is, $\Im h^*(z)=0$ for some
$z\in D$. Now it follows from the maximum principle that $h^*\equiv
0$ in $D$ and therefore in $M$. But then, as mentioned before, we
must have $\3{\6F}\cap (\{z\}\times\CC^N) = \{(z,p'_0)\}$ for all
$z\in M$, implying that $\2F$ and hence $F$ are constant. This
completes the proof. \epf

\medskip

\begin{rem} (a).
Assume that there exists $\eps>0$, such that for any $p\in M'$ and
$z$ in the ball $B_\eps(p)$, it holds that $M'\cap Q'_{z}\cap
B_\eps(p)=\{z\}$, for instance, if $M'$ is a strongly pseudoconvex
hypersurface. Then if $F_1(z)$ and $F_2(z)$  in
Proposition~\ref{KN-03}, are sufficiently close for some $z\in M\sm
{\cal S}$, it follows that $F_1 \equiv F_2$. In particular, $F$
cannot be extended as correspondence with a non-empty
(non-blowing-up) branch locus intersecting $M$.

(b). We also mention a paper by Shafikov in [Sha] where more
detailed studies in the equidimensional case ($N=n)$ were
addressed.

\end{rem}

\section{Hypersurfaces of positive Levi signature}\label{mixed}

 Fix
two integers $n,\ell$ with $1<\ell\le n/2$. For any $\epsilon$,
define
$${M_\epsilon}:=
\left\{[z_0,\cdots,z_{n+1}]\in \mathbf{CP}^{n+1}:
 |z|^2 \left(-\sum_{j=0}^{\ell} |z_j|^2 + \sum_{j=\ell+1}^{n+1}|z_j|^2\right)
 +\epsilon\left(|z_1|^4-|z_{n+1}|^4\right)
=0 \right\}.$$ Here $|z|^2=\sum_{j=0}^{n+1}|z_j|^2$ as usual. For
$\epsilon =0$, ${M_\epsilon}$ reduces to the generalized
sphere with signature $\ell$, which is the boundary of the
generalized ball
$$\mathbf {B}^{n+1}_\ell:=
\left\{\{[z_0,\cdots,z_{n+1}]\in \mathbf{CP}^{n+1}: -
\sum_{j=0}^{\ell}|z_j|^2 + \sum_{j=\ell+1}^{n+1}|z_j|^2 <0
\right\}.$$ The boundary $\p{\mathbf{B}_\ell^{n+1}}$  is locally
holomorphically equivalent to the hyperquadric $\mathbf
{H}^{n+1}_\ell\subset \mathbf{C}^{n+1}$ of signature $\ell$ defined
by $\Im{z_{n+1}}=-\sum_{j=1}^{\ell}|z_j|^2+
\sum_{j=\ell+1}^{n+1}|z_j|^2,$ where  $(z_1,\cdots, z_{n+1})$ is the
coordinates of $\mathbf{C}^{n+1}$.

For $0<\epsilon<<1$, ${M_\epsilon}$ is a compact smooth
real-algebraic hypersurface with Levi form nondegenerate of  the same signature
$\ell$. We now state our next theorem:

\bt\Label{2.1} There is an $\epsilon_0>0$ such that for
$0<\epsilon<\epsilon_0$, the following holds:
 ${M_\epsilon}$ is a
smooth real-algebraic hypersurface in $\mathbf{CP}^{n+1}$ with
nondegenerate Levi form of signature $\ell$ at every point. Moreover
for any point $p\in M_\epsilon$ and a holomorphic map $F$ from a
neighborhood $U$ of $p$ in $\mathbf {C}^{n+1}$ into $\mathbf
{C}^{N+1}$ sending $M_\epsilon$ into higher dimensional hyperquadric
$\mathbf {H}_\ell^{N+1}$ of the same signature, it follows that $F$
must be totally degenerate in the sense that $F({U})\subset \mathbf
{H}_\ell^{n+1}$. In particular, there does not exist  any
holomorphic embedding from any open piece of ${M_\e}$ into $ \mathbf
{H}_\ell^{N+1}$. \et

%

\bigskip
There are two main ingredients in our proof: the Chern-Moser-Weyl
tensor and an algebraicity  theorem of the first author in [Hu2]. We
first recall the related concept for the Chern-Moser-Weyl tensor.
For a more detailed account on this matter, the reader is referred
to [CM] and [HZ].

We use $(z,w)\in \mathbb{C}^{n}\ \times\mathbb{C}$ for the
coordinates of $\mathbb{C}^{n+1}$. We always assume that $n\ge 2$.
Let $M$ be a smooth real hypersurface.  $M$ is said to be  Levi
non-degenerate at $p\in M$ with signature $\ell\le n/2$ if there is
a local holomorphic change of coordinates, that maps $p$ to the
origin, such that in the new coordinates, $M$ is defined near $0$ by
an equation of the form:
\begin{equation}
r=v-|z|^2_\ell+o(|z|^2+|u|)=0. \label{eqn: 001}
\end{equation}
Here, we write $u=\Re w, v=\Im w$ and
$$\langle a,\bar b\rangle_\ell
=-\sum_{j=0}^\ell a_j \bar b_j+\sum_{j=\ell+1}^n a_j \bar b_j, \quad
|z|_\ell^2=\langle z,\bar z\rangle_\ell.$$
When $\ell=0$, we regard $\sum_{j\le
\ell} a_j=0$.

Assume that $M$ is Levi non-degenerate with the same signature
$\ell$ at any point with $\ell\le n/2$.  A contact form $\theta$
over $M$ is said to be appropriate if the Levi form $L_{\theta|_p}$
associated with $\theta$ at any  point $p\in M$ has $\ell$ negative
eigenvalues and $n-\ell$ positive eigenvalues. Let $\theta$ be an
appropriate contact form over $M$. Then from the Chern-Moser Theory,
there is a unique 4th order curvature tensor $\mathcal{S}_{\theta}$
associated with $\theta$ [CM], which we call the Chern-Moser-Weyl
tensor with respect to the contact form $\theta$ along $M$.
$\mathcal S_{\theta}$ can be regarded as a section in
$$T^{*(1,0)}M\otimes T^{*(0,1)}M\otimes T^{*(1,0)}M\otimes
T^{*(0,1)}M.$$
We write $\s S_{\theta|_p}$ for the restriction of
${\s S_{\theta}}$ at $p\in M$. For another appropriate contact for
$\wt{\theta}$, we have $\wt{\theta}=k\theta$ with $k\not = 0$.
Notice that $k>0$ when $\ell\not = n/2$. Then ${\s
S_{\wt\theta}}=k{\s S_{\theta}}$, i.e. the Chern-Moser-Weyl tensor
at a point $p\in M$ can be invariantly seen as multilinear map
$$\s S\colon T^{(1,0)}_pM\times T^{(0,1)}M_p\times T^{(1,0)}M_p\times
T^{(0,1)}_pM\to \C\otimes T_pM/ (T^{(1,0)}_pM+T^{(0,1)}_pM).$$

 The Chern-Moser-Weyl tensor has
particularly simple expression in the normal coordinates, which we
give as follws:
By the Chern-Moser normal form theory [CM], there is  a holomorphic
coordinates system in which $M$ is defined near $0$ by an equation
of the following form (see [(6.25), (6.30), CM]):
\begin{equation}
r=v-|z|_\ell^2+\frac{1}{4}s
(z,\bar{z})+o(|z|^4)=v-|z|_\ell^2+\frac{1}{4}\sum s_{\alpha\bar
{\beta}\gamma \bar{\delta}}z_{\alpha} {\bar z_\beta}z_{\gamma}{\bar
z_\delta}+o(|z|^4)=0, \label{eqn:002}
\end{equation}
where $s$ satisfies the trace condition
\begin{equation*}
\triangle_\ell s(z,\bar z)\equiv 0,
\end{equation*}
with $\triangle_\ell:=-\sum_{j\le \ell} \frac{\p^2}{\p z_j\p\bar
z_j} +\sum_{j=\ell+1}^n \frac{\p^2}{\p z_j\p\bar z_j}$. Here
$s(z,\-{z})=\sum s_{\alpha\bar \beta \gamma\bar{\delta}}z_{\alpha}
{\bar z_\beta}z_{\gamma}{\bar z_\delta}$, $\theta=i\p r,\
s_{\alpha\bar {\beta}\gamma \bar{\delta}}= s_{\gamma \bar
{\beta}\alpha\bar{\delta}}= s_{\gamma \bar{\delta}\alpha\bar
{\beta}},\ \overline {s_{\alpha\bar {\beta}\gamma
\bar{\delta}}}=s_{\beta\bar{\alpha}\delta\bar{\gamma}}$
 and
 the trace condition is equivalent to
  $\sum_{\alpha, \beta=1}^n s_{\alpha\bar {\beta}\gamma
\bar{\delta}}g^{\bar \beta \alpha}=0$ where $g^{\bar \beta
\alpha}=0$ for $\beta\neq\alpha$, $g^{\bar\beta\beta}=1$ for
$\beta>\ell, g^{\bar\beta\beta}=-1$ for $\beta\leq \ell$. Then
\begin{equation*}
s_{\alpha\bar {\beta}\gamma \bar{\delta}}=\s
S_{\theta|_0}(\frac{\p}{\p z_{\alpha}}|_0, {\frac{\p}{\p {\bar
z_{\beta}}}|_0}, \frac{\p}{\p z_{\gamma}}|_0, {\frac{\p}{\p  {\bar
z_{\delta}}}|_0}).
\end{equation*}
We also write $s_{\theta|_0}(z,\bar z)$ for $s(z,\-{z})$.
Consider the Levi null-cone
\begin{equation*}
\mathcal{C}_\ell=\{z\in{\CC}^n: |z|_\ell=0\}.
\end{equation*}
Then $\mathcal C_\ell$ is a real-algebraic variety of real
codimension $1$ in ${\CC}^n$ for $\ell>0$ with the only singularity at $0$.  For
each $p\in M$, write
$$\mathcal C_\ell T_p^{(1,0)}M=\{v_p\in
T_p^{(1,0)}M:\  \langle d\theta_p, v_p\wedge{\bar v_p} \rangle =0\}.$$
   $\mathcal
C_lT_p^{(1,0)}M$ is independent of the choice of $\theta$. Let $F$
be a CR diffeomorphism from $M$ to $M'$. We also have
$F_{*}(\mathcal C_\ell T_p^{(1,0)}M)=C_\ell T_{F(p)}^{(1,0)}M'$.
Write $\mathcal C_\ell T^{(1,0)}M=\coprod_{p\in M} \mathcal C_\ell
T_p^{(1,0)}M$ with the natural projection $\pi$ to $M$. We say that
$X$ is a smooth section of $\s C_\ell T^{(1,0)}M$ if $X$ is a smooth
vector field of type $(1,0)$ along $M$ such that $X|_p\in \s C_\ell
T_p^{(1,0)}M$ for each $p\in M$.

We say that the Chern-Moser-Weyl curvature tensor $\s S_{\theta}$ is
{\em pseudo positive semi-definite} (resp.\ {\em pseudo negative
semi-definite}) at $p\in M$ if $\s S_{\theta|_p}(X,\- X,X,\- X)\geq
0$  (resp.\ $\s S_{\theta|_p}(X,\-{X}, X, \- X)\leq 0$) for all
$X\in \s C_\ell T_p^{(1,0)}M$). We say that $\s S_{\theta}$ is {\em
pseudo positive definite} (resp.\ {\em pseudo negative definite}) at
$p\in M$ if $\s S_{\theta|_p}(X,\- X,X, \- X)> 0$  (resp.\ $\s
S_{\theta|_p}(X,\- X, X, \- X) < 0$)  for all $X\in \s C_\ell
T^{(1,0)}_pM \setminus{0}$). We use the terminology {\em pseudo
semi-definite} to mean either pseudo positive semi-definite or
pseudo negative semi-definite.

The following  will be used later:

\bt[\cite{hz1}, Corollary 3.3] \Label{2.3} Let $M\subset
{\CC}^{n+1}$ be a Levi non-degenerate real hypersurface of signature
$\ell$. Suppose that $F$ is a holomorphic mapping defined in a
(connected) open neighborhood $U$ of $M$ in ${\mathbf C}^{n+1}$ that
sends  $M$ into ${\mathbf H}_\ell ^{N+1}\subset {\CC}^{N+1}$. Assume
that $F(U)\not \subset {\mathbf H}_\ell^{N+1}$. Then when
$\ell<\frac{n}{2}$, the Chern-Moser-Weyl curvature tensor with
respect to any appropriate contact form $\theta$ is pseudo negative
semi-definite. When $\ell=\frac{n}{2}$, along any contact form
$\theta$, $\s S_{\theta}$ is pseudo semi-definite. \et

\bigskip
We next state the following algebraicity theorem:

\bt[\cite{Hu2}, Corollary  in $\S 2.3.5$] \Label{2.10}  Let
$M_1\subset {\mathbf C}^n$ and $M_2\subset {\mathbf C}^N$ with $N\ge
n\ge 2$ be two Levi non-degenerate real-algebraic hypersurfaces.
Let $p\in M_1$ and $U_p$ be a small
connected open neighborhood of $p$ in ${\mathbf C}^n$ and $F$ be a
holomorphic map from $U_p$ into ${\mathbf C}^N$ such that $F(U_p\cap
M_1)\subset M_2$ and $F(U_p)\not \subset M_2$. Suppose that $M_1$
and $M_2$ have the same signature $\ell$ at $p$ and $F(p)$,
respectively. Then $F$ is algebraic in the sense that each component
of $F$ satisfies a nontrivial holomorphic polynomial equation. \et

 It was first
proved in [Hu1] when $\ell=0$, namely, the strongly pseudo-convex
case. The case with $\ell>0$ was  done in the first author's thesis
[$\S 2.3.5$, Hu2]. It also follows from  a more general algebraicity
theorem of the second author in [Corollary1.6, Z1].

 The proof of the above theorem follows from the same  proof as in the
 signature zero case [Hu1], except in the $\ell>0$ case, we have now the Hopf lemma property as
 part of the assumption and that
the  proof of [Lemma 2.8, Hu1] (or [Lemma 2.8, Hu2]) needs to be
replaced by the following simple linear algebra lemma:


\bigskip
\bl[Lemma 2.8$'$, Hu2]
 Assume that $V$ is a smooth complex-analytic
hypersurface in a neighborhood of $0$ in ${\bf C}^{n+1}$. Assume
that $M'$ is a Levi non-degenerate hypersurface of signature
$\ell>0$ at $0$ and $T_0^{(1,0)}M'\not =T^{(1,0)}_0V$. Assume that
$M'\cap V$ contains a Levi non-degenerate CR submanifold of
hypersurface type with signature $\ell$ through $0$. Then $M'\cap V$
is a Levi non-degenerate hypersurface  of signature $\ell$ in $V$
near $0$.
\el
\bigskip


\bpf[Proof of Theorem \ref{2.1}]  We first note that for small
$\e$, $M_\e$ is a small perturbation of the compact quadric
$M_0$ of signature $\ell$. Hence there exists a positive
$0<\epsilon_0$ such that whenever $0<\e<\e_0$, ${M_\e}$ is
everywhere Levi non-degenerate with the same signature $\ell$.

Now, we compute the Chern-Moser-Weyl tensor of $M_\e$ at the point
$$P_0:=[\xi^0_0,\cdots,\xi^0_{n+1}], \quad \xi^0_j=0 \text{ for } j\not = 0,\ell+1, \quad \xi^0_0=1,\quad \xi_{\ell +1}^0=1,$$
 and consider the coordinates
$$\xi_0=1,\quad \xi_j=\frac{\eta_j}{1+\sigma}, \quad j=1,\cdots, \ell, \quad \xi_{\ell+1}=\frac{1-\sigma}{1+\sigma},
\quad \xi_{j+1}=\frac{\eta_j}{1+\sigma},\quad j=\ell+1 ,\cdots, n.$$
Then in the $(\eta,\sigma)$-coordinates, $P_0$ becomes the origin
and $M_\e$ is defined near the origin by an equation in the form:

\begin{equation}\label{22}
\rho=-4\Re{\sigma}-\sum_{j=1}^{\ell}|\eta_j|^2+\sum_{j=\ell+1}^{n}|\eta_j|^2
+{a}(|\eta_1|^4-|\eta_{n}|^4)+o(|\eta|^4)=0,
\end{equation}
for some $a>0$. Now, let $Q(\eta,\-\eta)=-a(|\eta_1|^4-|\eta_n|^4)$
and make a standard $\ell$-harmonic decomposition [SW]:

\begin{equation}\label {33}
Q(\eta,\-\eta)=N^{(2,2)}(\eta,\-\eta)+
A^{(1,1)}(\eta,\-\eta)|\eta|^2_{\ell}.
\end{equation}
Here $N^{(2,2)}(\eta,\eta)$ is a $(2,2)$-homogeneous polynomial in
$(\eta,\-{\eta})$ such that $\Delta_\ell N^{(2,2)}(\eta,\-\eta)=0$
with $\Delta_\ell$ as before. Now $N^{(2,2)}$ is the
Chern-Moser-Weyl tensor of $M_\e$ at $0$ (with respect to an obvious
contact form) with $N^{(2,2)}(\eta,\-\eta)=Q(\eta,\-\eta)$ for any
$\eta\in{\mathcal C}T^{(1,0)}_0 M_e$. Now the value of the
Chern-Moser-Weyl tension has negative  and positive value at
$X_1=\frac{\p }{\p \eta_1}+\frac{\p }{\p \eta_{\ell+1}}|_0$ and
$X_2=\frac{\p }{\p \eta_{2}}+\frac{\p }{\p \eta_{n}}|_0$,
respectively. If $\ell>1$, then  both $X_1$ and $X_2$ are in
${\mathcal C}T^{(1,0)}_0 M_e$. We see that the Chern-Moser-Weyl
tensor can not be semi-definite near the origin in such a coordinate
system.


Next, suppose an open piece $U$ of ${M_\e}$ can be holomorphically
and transversally embedded into the ${\mathbf H}_\ell^{N+1}$ for
$N>n$ by $F$. Then by the algebraicity result in Theorem \ref{2.10},
$F$ is algebraic. Since the branching points of $F$ and the points
where $F$ is not defined (poles or points of indeterminancy of $F$)
are contained in a complex-algebraic variety of codimension at most
one, $F$ extends holomorphically along a smooth curve  $\gamma$
starting from some point in $U$ and ending up at some point $p^*
(\approx 0)\in M_\e$ in the $(\eta,\sigma)$-space where the
Chern-Moser-Weyl tensor of $M_\e$ is not pseudo-semi-definite. By
the uniqueness of real-analytic functions, the extension of $F$ must
also map an open piece of $p^*$ into  $ {\bf H}^{N+1}_\ell$. The
extension is not totally degenerate. By Theorem \ref{2.3}, we get a
contradiction. \epf

\bigskip

\section{Open problems}
We mention here the following questions that still
seem to be  open.

\medskip
{\bf Question 1:} {\em Is there any example of a compact strongly pseudoconvex
real-algebraic hypersurface in $\C^n$ that is not holomorphically embeddable into a sphere of any dimension?}
\medskip

In fact, all known examples of hypersurfaces that are not
embeddable into spheres are also not embeddable into strongly pseudoconvex real-algebraic hypersurfaces.
It remains unknown whether these two classes are different,
more precisely:

\medskip
{\bf Question 2:} {\em Is there any example of a (not necessarily compact) strongly pseudoconvex
real-algebraic hypersurface in $\C^n$ that is holomorphically embeddable
into a compact strongly pseudoconvex real-algebraic hypersurface
but is not holomorphically embeddable into a sphere of any dimension?}
\bigskip

{\bf Acknowledgement}:
This work started   from the authors'
 examples (those given in $\S
2$)  observed during an International Conference on Geometric
Analysis in Several Complex Variables and PDEs held at Serra Negra,
Brazil in August 2011. Both authors thank the organizers of the
conference for the invitation and for providing a stimulating
academical environment. The authors also would like to thank the
referee for a  careful reading and many comments which helped
greatly to improve the exposition of the paper.

\end{document}